\documentclass[a4paper,12pt,british,english]{article}
\usepackage[a4paper, margin=2cm]{geometry}
\usepackage[T1]{fontenc}
\usepackage[latin9]{inputenc}
\usepackage{color}
\usepackage{amsmath}
\usepackage{amssymb}
\usepackage{esint}
\PassOptionsToPackage{normalem}{ulem}
\usepackage{ulem}
\usepackage{enumerate}

\usepackage{babel}

\newtheorem{theo}{Theorem}

\newtheorem{prop}{Proposition}  
\newtheorem{coro}{Corollary}

\newtheorem{rem}{Remark}
\newtheorem{lema}{Lemma}
\newtheorem{fact}{Fact}
\newtheorem{defi}{Definition}

\newenvironment{dwd}{\par\noindent{\bf Proof.}}{\par\rightline{$\blacksquare$}}

 \global\long\def\sbr#1{\left[ #1\right] }
 \global\long\def\cbr#1{\left\{  #1\right\}  }
 \global\long\def\rbr#1{\left(#1\right)}

 \global\long\def\R{\mathbb{R}}

 \global\long\def\TDD#1{}
 \global\long\def\dd#1{\textnormal{d}#1}

 \global\long\def\TTV#1#2#3{\text{TV}^{#3}\!\rbr{#1,#2}}

 \global\long\def\ra{\rightarrow}
 
 \global\long\def\ns{\infty}

\global\long\def\Varnorm#1#2{\left\Vert {#1} \right\Vert_{{#2}-\text{var},\left[a;b\right]}}
\global\long\def\Varnormthm#1#2{\left\Vert {#1} \right\Vert_{{#2}-\emph{var},\left[a;b\right]}}
\global\long\def\Varfullnorm#1#2{\left\Vert {#1} \right\Vert_{\text{var},{#2},\left[a;b\right]}}
\global\long\def\Varfullnormthm#1#2{\left\Vert {#1} \right\Vert_{\emph{var},{#2},\left[a;b\right]}}
\global\long\def\Oscnorm#1{\left\Vert {#1} \right\Vert_{\text{osc},\left[a;b\right]}}

\begin{document}

\title{A new theorem on the existence of the Riemann-Stieltjes integral and an improved version of the Lo\'{e}ve-Young inequality}

\author{Rafa{\l{}} M. \L{}ochowski \\ \\
{\small Department of Mathematics and Mathematical Economics, Warsaw School of Economics} \\
{\small ul. Madali\'{n}skiego 6/8, 02-513 Warszawa, Poland}\\
{\small E-mail : rlocho314@gmail.com}}

\date{}

\maketitle

\begin{abstract}
Using the notion of truncated variation we obtain a new theorem on the existence and estimation of the Riemann-Stieltjes integral. As a special case of this theorem we obtain an improved version of the Lo\'{e}ve-Young inequality for the Riemann-Stieltjes integrals driven by irregular signals. Using this result we strenghten some results of Terry Lyons on the existence of solutions of integral equations driven by moderately irregular signals. 

\end{abstract}

\section{Introduction}

The purpose of this paper is to investigate the top-down structure of the Riemann-Stieltjes integral and to state some general condition guaranteeing the existence of this integral, expressed in terms of the functional called {\em truncated variation}. 

The simplest (and rather not surprising) case where the Riemann-Stieltjes integral  $\int_{a}^{b}f\mathrm{d}g$ (RSI in short) exists, is the situation when the integrand and the integrator have no common points of discontinuity, the former is bounded and the latter has finite total variation. We will prove a general theorem (Theorem \ref{main}) encompassing this situation as well as more interesting case when the RSI exists, namely when the integrand
and the integrator have possibly unbounded variation, but they have finite
$p$-variation and $q$-variation respectively, with $p>1,$ $q>1$
and $p^{-1}+q^{-1}>1.$ The latter result is due to Young (\cite[p. 264, Theorem on Stieltjes integrability]{Young:1936}). For $f:\left[a;b\right]\rightarrow\R$
and $p>0,$ the $p$-variation, which we will denote by $V^{p}\left(f;\left[a;b\right]\right),$
is defined as 
\[
V^{p}\left(f,\left[a;b\right]\right)=\sup_{n}\sup_{a\leq t_{1}<t_{2}<\cdots<t_{n}\leq b}\sum_{i=1}^{n-1}\left|f\left(t_{i}\right)-f\left(t_{i-1}\right)\right|^{p}.
\]
The aforementioned Theorem \ref{main} provides also an upper bound for the difference 
\[
\left| \int_{a}^{b}f\mathrm{d}g - f(a)\sbr{g(b)-g(a)} \right|.
\]
As a special case of this bound we will obtain the following, improved version of the classical Lo\'{e}ve-Young inequality:
\begin{align} \label{LYimproved}
\left|\int_{a}^{b}f\mathrm{d}g-f\left(a\right)\left[g\left(b\right)-g\left(a\right)\right]\right| 
  \leq C_{p,q}\rbr{V^{p}\left(f,\left[a;b\right]\right)}^{1-1/q}\left\Vert f\right\Vert _{\text{osc},\left[a;b\right]}^{1+p/q-p} \rbr{V^{q}\left(g,\left[a;b\right]\right)}^{1/q},
\end{align}
where $C_{p,q}$ is some constant depending on $p$ and $q$ only and $\left\Vert f\right\Vert_{\text{osc},\left[a;b\right]} := \sup_{a \leq s < t \leq b} \left| f(t) - f(s)\right|.$
The original Lo\'{e}ve-Young estimate, published in 1936 in \cite{Young:1936}, reads as:
\begin{equation*}
\left| \int_{a}^{b}f\mathrm{d}g - f(a)\sbr{g(b)-g(a)} \right| \leq \zeta\rbr{p^{-1}+q^{-1}} \rbr{V^{p}\left(f,\left[a;b\right]\right)}^{1/p} \rbr{V^{q}\left(g,\left[a;b\right]\right)}^{1/q},
\end{equation*} 
where $\zeta \rbr{r} = \sum_{k=1}^{+\infty} k^{-r}$ is the famous Riemann zeta function.
We say that our bound is improved version of the Lo\'{e}ve-Young inequality since the constant $C_{p,q}$ (although possibly greater than $\zeta\rbr{p^{-1}+q^{-1}}$) is irrelevant in applications while 
the fact that the term $\rbr{V^{p}\left(f,\left[a;b\right]\right)}^{1/p}$ in the original Lo\'{e}ve-Young estimate is replaced in our estimate by $\rbr{V^{p}\left(f,\left[a;b\right]\right)}^{1-1/q}\left\Vert f\right\Vert _{\text{osc},\left[a;b\right]}^{1+p/q-p}$ (notice that $1/p >1-1/q $ and always $\rbr{V^{p}\left(f,\left[a;b\right]\right)}^{1/p-\rbr{1-1/q}} \geq \left\Vert f\right\Vert _{\text{osc},\left[a;b\right]}^{1+p/q-p} $) makes it possible to obtain stronger results. For example, Proposition \ref{prop:int_equation} obtained with the help of (\ref{LYimproved}) is a genuine improvement upon earlier known results of this type. 

Let us comment shortly how the results on the existence of the RSI were obtained so far. The original Young's proof utlilized elementary but clever induction
argument for finite sequences. Another proof of the Young theorem may
be found in \cite[Chapt. 6]{Friz:2010fk}, where integral estimates
based on a control function and the Lóeve-Young inequality are used. This
approach is further applied in the rough-path theory setting. Further generalisations of Young's theorem are possible, with $p$-variation replaced by more general $\varphi$-variation: 
\[
V^{\varphi}\left(f,\left[a;b\right]\right)=\sup_{n}\sup_{a\leq t_{1}<t_{2}<\cdots<t_{n}\leq b}\sum_{i=1}^{n-1} \varphi \rbr{\left|f\left(t_{i}\right)-f\left(t_{i-1}\right)\right|},
\]
where $\varphi:[0;+\ns) \ra [0;+\ns)$ is a Young function, i.e. convex, strictly increasing function starting from $0$ (see for example \cite{Young:1938}, \cite{Dyackov:1988} and for a survey about another results of this type see the recent books \cite[Chapt. 3]{NorvaisaConcrete:2010},  \cite[Sect. 4.4]{Banas_et_al:2013}).

However, as far as we know, Theorem \ref{main} is a new result on the existence of the RSI. The proof of Theorem \ref{main} utilizes simple properties of the truncated variation and multiple application of the summation by parts.
Similarly, no version of  the Lo\'{e}ve-Young inequality as estimate (\ref{LYimproved}), as far as we know, has appeared so far (see detailed historical notes on the the Lo\'{e}ve-Young inequality in \cite[pp. 212-214]{NorvaisaConcrete:2010}). We conjecture that using Theorem \ref{main} one may also obtain a variation of the  Lo\'{e}ve-Young inequality for $\varphi$-variation (see \cite[Theorem 3.89, Corollary 3.90]{NorvaisaConcrete:2010} or \cite[Theorem 4.40]{Banas_et_al:2013}). We intend to deal with this conjecture in the future. 

After having obtained these results we were able, following Lyons \cite{Lyons:1994}, and Lyons, Caruana and L\'{e}vy \cite{LyonsCaruana:2007}, to solve few types of integral equations driven by moderately irregular signals. By moderately irregular signals we mean continuous signals with finite $p$-variation, where $p\in (1;2).$ It is well known that for higher degrees of irregularity, corresponding to $p\geq 2,$ one needs, constructing approximations of integral equations, to consider terms of a new type (like L\'{e}vy's area). We believe that the tuncated variation approach for such paths is also possible and this will be a topic of our further research. 

Let us comment shortly on the organisation of the paper. In the next section we prove a general theorem on the existence of the Riemann-Stieltjes integral, expressed in terms of the truncated variation functionals and derive from it the stronger version of the Lo\'{e}ve-Young inequality. Next, in Section 3, we deal with the applications of this result to few types of integral equations driven by moderately irregular signals.
 
\section{A theorem on the existence of the Riemann-Stieltjes integral} 
\label{RSI_existence}

In this section we will prove a general theorem on the existence of the RSI $\int_{a}^{b}f\mathrm{d}g$ formulated in terms of the truncated variation.
We will assume that both - integrand $f:[a;b]\ra \R$ and integrator $g:[a;b]\ra \R$ are regulated functions. Let us state the definition of the truncated variation and recall the definition of a regulated function. 

For $f:[a;b]\ra \R$  its {\em truncated variation} with the truncation parameter $\delta \geq 0$ will be denoted by $\TTV f{[a;b]}{\delta}.$ It may be simply defined as the greatest lower bound for the total variation of any function $g:[a;b]\ra \R,$ uniformly approximating $f$ with the accuracy $\delta/2,$
\[
\TTV f{[a;b]}{\delta} := \inf\cbr{\TTV {g}{[a;b]}{}: \|f - g \|_{\ns,[a;b]} \leq \delta/2}.
\]
$\|f - g \|_{\ns, [a;b]}$ denotes here $\sup_{a\leq t \leq b} \left|f(t) - g(t)\right|$ and the total variation $\TTV {g}{[a;b]}{}$ is defined as 
\[\TTV {g}{[a;b]}{} := \sup_{n}\sup_{a\leq t_{1}<t_{2}<\cdots<t_{n}\leq b}\sum_{i=1}^{n-1}\left|g\left(t_{i}\right)-g\left(t_{i-1}\right)\right|.  \]
It appears that the truncated variation $\TTV f{[a;b]}{\delta}$ is finite for any $\delta>0$ iff $f$ is {\em regulated} (cf. \cite[Fact 2.2]{LochowskiColloquium:2013}) and then for any $\delta >0$  the following equality holds
\begin{equation} \label{TV_def1}
\TTV f{[a;b]}{\delta} = \sup_{n}\sup_{a\leq t_{1}<t_{2}<\cdots<t_{n}\leq b}\sum_{i=1}^{n-1}\max\left\{ \left|f\left(t_{i}\right)-f\left(t_{i-1}\right)\right|-\delta,0\right\}
\end{equation}
(cf. \cite[Theorem 4]{LochowskiGhomrasniMMAS:2014}). 

Let us recall that a function  $h:[a;b]\ra \R$ is {\em regulated} if there exist one-sided finite limits
$\lim_{t\rightarrow a+}h\left(t\right)$ and $\lim_{t\rightarrow b-}h\left(t\right),$
and for any $t\in\left(a;b\right)$ and there exist one-sided finite limits $\lim_{t\rightarrow x-}h\left(t\right)$
and $\lim_{t\rightarrow x+}h\left(t\right).$

We will also need the following
result (cf. \cite[Theorem 4]{LochowskiGhomrasniMMAS:2014}): for any regulated
function $f:\left[a;b\right]\rightarrow\R$ and $\delta>0$ there
exists a regulated function $f^{\delta}:\left[a;b\right]\rightarrow\R$
such that $\left\Vert f-f^{\delta}\right\Vert _{\infty,\left[a;b\right]}\leq\delta/2$
and 
\[
\TTV{f^{\delta}}{\left[a;b\right]}0=\TTV f{\left[a;b\right]}{\delta}.
\]
From formula (\ref{TV_def1}), it directly follows that the truncated variation
is a superadditive functional of the interval, i.e. for any $d\in\left(a;b\right)$
\[
\TTV f{\left[a;b\right]}{\delta}\geq\TTV f{\left[a;d\right]}{\delta}+\TTV f{\left[d;b\right]}{\delta}.
\]
Moreover, we also have the following easy estimate of the truncated
variation of a function $f$ perturbed by some other function $h:$
\begin{equation}
\TTV{f+h}{\left[a;b\right]}{\delta}\leq\TTV f{\left[a;b\right]}{\delta}+\TTV h{\left[a;b\right]}0,\label{eq:TV_variation}
\end{equation}
which stems directly from the inequality: for $a\leq s<t\leq b,$
\begin{align*} & \max\left\{ \left|f\left(t\right)+h\left(t\right)-\left\{ f\left(s\right)+h\left(s\right)\right\} \right|-\delta, 0\right\}  \\ & \quad \quad \quad \quad \quad \quad \quad \quad \leq \max\left\{ \left|f\left(t\right)-f\left(s\right)\right|-\delta,0\right\} +\left|h\left(t\right)-h\left(s\right)\right|.
\end{align*}

\begin{theo} \label{main} Let $f,g:\left[a;b\right]\rightarrow\R$
be two regulated functions which have no common points of discontinuity.
Let $\eta_{0}\geq\eta_{1}\geq\ldots$ and $\theta_{0}\geq\theta_{1}\geq\ldots$
be two sequences of non-negative numbers, such that $\eta_{k}\downarrow0,$
$\theta_{k}\downarrow0$ as $k\rightarrow+\infty.$ Define $\eta_{-1}:=\sup_{a\leq t\leq b}\left|f\left(t\right)-f\left(a\right)\right|$
and 
\[
S:=\sum_{k=0}^{+\infty}2^{k}\eta_{k-1}\cdot\TTV g{\left[a;b\right]}{\theta_{k}}+\sum_{k=0}^{\infty}2^{k}\theta_{k}\cdot\TTV f{\left[a;b\right]}{\eta_{k}}.
\]
If $S<+\infty$ then the Riemann-Stieltjes integral 
$
\int_{a}^{b}f\mathrm{d}g
$
exists and one has the following estimate 
\begin{equation}
\left|\int_{a}^{b}f\mathrm{d}g-f\left(a\right)\left[g\left(b\right)-g\left(a\right)\right]\right|\leq S.\label{eq:estimate_integral}
\end{equation}
\end{theo}

\begin{rem} The assumption that $f$ and $g$ has no common points
of discontinuity is necessary for the existence of the RSI $\int_{a}^{b}f\mathrm{d}g.$
When a more general integrals are considered (e.g. the Moore-Pollard
integral, c.f. \cite[p. 263]{Young:1936}), we may weaken this assumption
and assume that $f$ and $g$ have no common one-sided discontinuities.
\end{rem} The proof of Theorem \ref{main} will be based on the following
lemma.

\begin{lema} \label{lema} Let $f,g:\left[a;b\right]\rightarrow\R$
be two regulated functions. Let $c=t_{0}<t_{1}<\ldots<t_{n}=d$ be
any partition of the interval $\left[c;d\right]\subset\left[a;b\right]$
and let $\xi_{0}=c$ and $\xi_{1},\ldots,\xi_{n}$ be such that $t_{i-1}\leq\xi_{i}\leq t_{i}$
for $i=1,2,\ldots,n.$ Then for $\delta_{-1}:=\sup_{c\leq t\leq d}\left|f\left(t\right)-f\left(c\right)\right|,$
 $\delta_{0}\geq\delta_{1}\geq\ldots\geq\delta_{r}\geq0$ and
$\varepsilon_{0}\geq\varepsilon_{1}\geq\ldots\geq\varepsilon_{r}\geq0$
the following estimate holds 
\begin{eqnarray*}
 &  & \left|\sum_{i=1}^{n}f\left(\xi_{i}\right)\left[g\left(t_{i}\right)-g\left(t_{i-1}\right)\right]-f\left(c\right)\left[g\left(d\right)-g\left(c\right)\right]\right|\\
 &  & \leq\sum_{k=0}^{r}2^{k}\delta_{k-1}\cdot\TTV g{\left[c;d\right]}{\varepsilon_{k}}+\sum_{k=0}^{r}2^{k}\varepsilon_{k}\cdot\TTV f{\left[c;d\right]}{\delta_{k}}+n\delta_{r}\varepsilon_{r}.
\end{eqnarray*}
\end{lema} 
\begin{dwd} Denote $\varepsilon=\varepsilon_{0},$ by
summation by parts, we have the following equality 
\begin{align}
 & \sum_{i=1}^{n}  f\left(\xi_{i}\right)\left[g\left(t_{i}\right)-g\left(t_{i-1}\right)\right]-f\left(c\right)\left[g\left(d\right)-g\left(c\right)\right]  \nonumber \\
 & =  \sum_{i=1}^{n}\left[f\left(\xi_{i}\right)-f\left(c\right)\right]\left[g^{\varepsilon}\left(t_{i}\right)-g^{\varepsilon}\left(t_{i-1}\right)\right]\nonumber \\
 &  \quad+\sum_{i=1}^{n}\left[f\left(\xi_{i}\right)-f\left(c\right)\right]\left[g\left(t_{i}\right)-g^{\varepsilon}\left(t_{i}\right)-\left\{ g\left(t_{i-1}\right)-g^{\varepsilon}\left(t_{i-1}\right)\right\} \right]\nonumber \\
 & =  \sum_{i=1}^{n}\left[f\left(\xi_{i}\right)-f\left(c\right)\right]\left[g^{\varepsilon}\left(t_{i}\right)-g^{\varepsilon}\left(t_{i-1}\right)\right]\nonumber \\
 & \quad +\sum_{i=1}^{n}\left[g\left(d\right)-g^{\varepsilon}\left(d\right)-\left\{ g\left(t_{i-1}\right)-g^{\varepsilon}\left(t_{i-1}\right)\right\} \right]\left[f\left(\xi_{i}\right)-f\left(\xi_{i-1}\right)\right],\label{eq:jeden}
\end{align}
where $g^{\varepsilon}:[c;d]\ra \R$ is regulated and such that 
\[
\left\Vert g-g^{\varepsilon}\right\Vert _{\infty,[c;d]}\leq\frac{1}{2}\varepsilon\mbox{ and }\TTV{g^{\varepsilon}}{\left[c;d\right]}0=\TTV g{\left[c;d\right]}{\varepsilon}.
\]
Similarly, for $\delta=\delta_{0}$ we may write 
\begin{align}
 & \sum_{i=1}^{n} \left[g\left(d\right)-g^{\varepsilon}\left(d\right)-\left\{ g\left(t_{i-1}\right)-g^{\varepsilon}\left(t_{i-1}\right)\right\} \right]\left[f\left(\xi_{i}\right)-f\left(\xi_{i-1}\right)\right]\nonumber \\
 & =  \sum_{i=1}^{n}\left[g\left(d\right)-g^{\varepsilon}\left(d\right)-\left\{ g\left(t_{i-1}\right)-g^{\varepsilon}\left(t_{i-1}\right)\right\} \right]\left[f^{\delta}\left(\xi_{i}\right)-f^{\delta}\left(\xi_{i-1}\right)\right]\label{eq:dwa}\\
 &  \quad +\sum_{i=1}^{n}\left[f\left(\xi_{i}\right)-f^{\delta}\left(\xi_{i}\right)-\left\{ f\left(c\right)-f^{\delta}\left(c\right)\right\} \right]\left[\left\{ g\left(t_{i}\right)-g^{\varepsilon}\left(t_{i}\right)\right\} -\left\{ g\left(t_{i-1}\right)-g^{\varepsilon}\left(t_{i-1}\right)\right\} \right],\nonumber 
\end{align}
where $f^{\delta}:[c;d]\ra \R$ is regulated and such that 
\[
\left\Vert f-f^{\delta}\right\Vert _{\infty,[c;d]}\leq\frac{1}{2}\delta\mbox{ and }\TTV{f^{\delta}}{\left[c;d\right]}0=\TTV f{\left[c;d\right]}{\delta}.
\]
Since $\TTV{g^{\varepsilon}}{\left[c;d\right]}0=\TTV g{\left[c;d\right]}{\varepsilon},$
$\TTV{f^{\delta}}{\left[c;b\right]}0=\TTV f{\left[c;d\right]}{\delta},$
$\left\Vert g-g^{\varepsilon}\right\Vert _{\infty,[c;d]}\leq\varepsilon/2$
and $\left\Vert f-f^{\delta}\right\Vert _{\infty,[c;d]}\leq\delta/2,$
from (\ref{eq:jeden}) and (\ref{eq:dwa}) we have the following estimate
\begin{eqnarray}
 &  & \left|\sum_{i=1}^{n}f\left(\xi_{i}\right)\left[g\left(t_{i}\right)-g\left(t_{i-1}\right)\right]-f\left(c\right)\left[g\left(d\right)-g\left(c\right)\right]\right|\nonumber \\
 &  & \leq\sup_{c\leq t\leq d}\left|f\left(t\right)-f\left(c\right)\right|\cdot\TTV g{\left[c;d\right]}{\varepsilon}+\varepsilon\cdot\TTV f{\left[c;d\right]}{\delta}+n\delta\varepsilon.\label{eq:jeden1}
\end{eqnarray}
Denote $g_{1}:=g-g^{\varepsilon},$ $f_{1}:=f-f^{\delta}$ on $[c;d].$ By (\ref{eq:jeden})
and (\ref{eq:dwa}), instead of the last summand $n\delta\varepsilon$
in (\ref{eq:jeden1}) we may write the estimate 
\begin{eqnarray}
 &  & \left|\sum_{i=1}^{n}\left[f\left(\xi_{i}\right)-f^{\delta}\left(\xi_{i}\right)-\left\{ f\left(c\right)-f^{\delta}\left(c\right)\right\} \right]\left[\left\{ g\left(t_{i}\right)-g^{\varepsilon}\left(t_{i}\right)\right\} -\left\{ g\left(t_{i-1}\right)-g^{\varepsilon}\left(t_{i-1}\right)\right\} \right]\right|\nonumber \\
 &  & =\left|\sum_{i=1}^{n}\left[f_{1}\left(\xi_{i}\right)-f_{1}\left(c\right)\right]\left[g_{1}\left(t_{i}\right)-g_{1}\left(t_{i-1}\right)\right]\right|\nonumber \\
 &  & \leq\sup_{c\leq t\leq d}\left|f_{1}\left(t\right)-f_{1}\left(c\right)\right|\cdot\TTV{g_{1}}{\left[c;d\right]}{\varepsilon_{1}}+\varepsilon_{1}\cdot\TTV{f_{1}}{\left[c;d\right]}{\delta_{1}}+n\delta_{1}\varepsilon_{1}\nonumber \\
 &  & \leq\delta\cdot\TTV{g_{1}}{\left[c;d\right]}{\varepsilon_{1}}+\varepsilon_{1}\cdot\TTV{f_{1}}{\left[c;d\right]}{\delta_{1}}+n\delta_{1}\varepsilon_{1},\label{eq:dwa1}
\end{eqnarray}
where the last but one inequlity in (\ref{eq:dwa1}) follows by the
same reasoning for $f_{1}$ and $g_{1}$ as inequality (\ref{eq:jeden})
for $f$ and $g.$ Repeating these arguments, by induction we get
\begin{eqnarray}
 &  & \left|\sum_{i=1}^{n}f\left(\xi_{i}\right)\left[g\left(t_{i}\right)-g\left(t_{i-1}\right)\right]-f\left(c\right)\left[g\left(d\right)-g\left(c\right)\right]\right|\nonumber \\
 &  & \leq\sum_{k=0}^{r}\delta_{k-1}\cdot\TTV{g_{k}}{\left[c;d\right]}{\varepsilon_{k}}+\sum_{k=0}^{r}\varepsilon_{k}\cdot\TTV{f_{k}}{\left[c;d\right]}{\delta_{k}}+n\delta_{r}\varepsilon_{r},\label{eq:induction}
\end{eqnarray}
where $\delta_{-1}:=\sup_{c\leq t\leq c}\left|f\left(t\right)-f\left(a\right)\right|,$
$g_{0}\equiv g,$ $f_{0}\equiv f$ and for $k=1,2,\ldots,r,$ $g_{k}:=g_{k-1}-g_{k-1}^{\varepsilon_{k-1}},$
$f_{k}:=f_{k-1}-f_{k-1}^{\delta_{k-1}}$ are defined similarly as
$g_{1}$ and $f_{1}.$\par Since $\varepsilon_{k}\leq\varepsilon_{k-1}$
for $k=1,2,\ldots,r,$ by (\ref{eq:TV_variation}) and the fact that
the function $\delta\mapsto\TTV h{\left[c;d\right]}{\delta}$ is non-increasing,
we estimate 
\begin{eqnarray*}
\TTV{g_{k}}{\left[c;d\right]}{\varepsilon_{k}} & = & \TTV{g_{k-1}-g_{k-1}^{\varepsilon_{k-1}}}{\left[c;d\right]}{\varepsilon_{k}}\\
 & \leq & \TTV{g_{k-1}}{\left[c;d\right]}{\varepsilon_{k}}+\TTV{g_{k-1}^{\varepsilon_{k-1}}}{\left[c;d\right]}0\\
 & = & \TTV{g_{k-1}}{\left[c;d\right]}{\varepsilon_{k}}+\TTV{g_{k-1}}{\left[c;d\right]}{\varepsilon_{k-1}}\\
 & \leq & 2\TTV{g_{k-1}}{\left[c;d\right]}{\varepsilon_{k}}.
\end{eqnarray*}
Hence, by recursion, for $k=1,2,\ldots,r,$ 
\[
\TTV{g_{k}}{\left[c;d\right]}{\varepsilon_{k}}\leq2^{k}\TTV g{\left[c;d\right]}{\varepsilon_{k}}.
\]
Similarly, for $k=1,2,\ldots,r,$ we have 
\[
\TTV{f_{k}}{\left[c;d\right]}{\delta_{k}}\leq2^{k}\TTV f{\left[c;d\right]}{\delta_{k}}.
\]
By (\ref{eq:induction}) and last two estimates we get the desired
estimate.
\end{dwd} 
\begin{rem} \label{symmetry}
Notice that starting in (\ref{eq:jeden}) from the summation by parts, then splitting the difference $f\left(\xi_{i}\right)-f\left(\xi_{i-1}\right):$
\begin{align*}
 & \sum_{i=1}^{n}  f\left(\xi_{i}\right)\left[g\left(t_{i}\right)-g\left(t_{i-1}\right)\right]-f\left(c\right)\left[g\left(d\right)-g\left(c\right)\right]  \\
 & =  \sum_{i=1}^{n}\left[g\left(d\right)-g\left(t_{i-1}\right)\right]
\left[f\left(\xi_{i}\right)-f\left(\xi_{i-1}\right)\right]
\\
 & =  \sum_{i=1}^{n}\left[g\left(d\right)-g\left(t_{i-1}\right)\right]\left[f^{\delta}\left(\xi_{i}\right)-f^{\delta}\left(\xi_{i-1}\right)\right] \\
& \quad + \sum_{i=1}^{n}\left[g\left(d\right)-g\left(t_{i-1}\right)\right]\left[f\left(\xi_{i}\right)-f^{\delta}\left(\xi_{i}\right)
-\cbr{f\left(\xi_{i-1}\right) - f^{\delta}\left(\xi_{i-1}\right) }\right] 
\end{align*}
and proceeding similarly as in the proof of Lemma \ref{lema} we get the symmetric estimate 
\begin{eqnarray}
 &  & \left|\sum_{i=1}^{n}f\left(\xi_{i}\right)\left[g\left(t_{i}\right)-g\left(t_{i-1}\right)\right]-f\left(c\right)\left[g\left(d\right)-g\left(c\right)\right]\right| \nonumber \\
 &  & \leq\sum_{k=0}^{r}2^{k}\varepsilon_{k-1}\cdot\TTV f{\left[c;d\right]}{\delta_{k}}+\sum_{k=0}^{r}2^{k}\delta_{k}\cdot\TTV g{\left[c;d\right]}{\varepsilon_{k}}+n\delta_{r}\varepsilon_{r}, \label{symmetry_left}
\end{eqnarray}
where $\varepsilon_{-1}=\sup_{c\leq t \leq d}\left|g(d)-g(t) \right|.$
\end{rem}
\begin{rem} \label{n=1}
Setting in Lemma \ref{lema}, $n=1$ for any $\xi \in [c;d]$ we get the estimate 
\begin{eqnarray}
 &  & \left|\left(f\left(\xi\right)-f\left(c\right)\right)\left[g\left(d\right)-g\left(c\right)\right]\right| \nonumber \\
 &  & \leq\sum_{k=0}^{r}2^{k}\delta_{k-1}\cdot\TTV f{\left[c;d\right]}{\varepsilon_{k}}+\sum_{k=0}^{r}2^{k}\varepsilon_{k}\cdot\TTV g{\left[c;d\right]}{\delta_{k}}+n\delta_{r}\varepsilon_{r}. \label{n=1_left}
\end{eqnarray}
and similarly, setting in Remark \ref{symmetry}, $n=1$ we get similar estimate, where the right side of (\ref{n=1_left}) is replaced by the right side of  (\ref{symmetry_left}).
\end{rem}
Now we proceed to the proof of Theorem \ref{main}.
\begin{dwd} It is enough to prove that for any two partitions 
\[
\pi=\left\{ a=a_{0}<a_{1}<\ldots<a_{l}=b\right\} ,
\]
\[
\rho=\left\{ a=b_{0}<b_{1}<\ldots<b_{m}=b\right\} 
\]
and $\nu_{i}\in\left[a_{i-1};a_{i}\right],$ $\xi_{j}\in\left[b_{j-1};b_{j}\right],$
$i=1,2,\ldots,l,$ $j=1,2,\ldots,m,$ the difference 
\[
\left|\sum_{i=1}^{l}f\left(\nu_{i}\right)\left[g\left(a_{i}\right)-g\left(a_{i-1}\right)\right]-\sum_{j=1}^{m}f\left(\xi_{j}\right)\left[g\left(b_{j}\right)-g\left(b_{j-1}\right)\right]\right|
\]
is as small as we please, provided that the meshes of the partitions
$\pi$ and $\rho$, defined as 
\[
\mbox{mesh}\left(\pi\right):=\max_{i=1,2,\ldots,l}\left(a_{i}-a_{i-1}\right),\text{ }\mbox{mesh}\left(\rho\right):=\max_{j=1,2,\ldots,m}\left(b_{j}-b_{j-1}\right)
\]
respectively, are sufficiently small. 

Define 
\[
\sigma=\pi\cup\rho=\left\{ a=s_{0}<s_{1}<\ldots<s_{n}=b\right\} 
\]
and for $i=1,2,\ldots,l$ consider 
\[
\left|f\left(\nu_{i}\right)\left[g\left(a_{i}\right)-g\left(a_{i-1}\right)\right]-\sum_{k:s_{k-1},s_{k}\in\left[a_{i-1};a_{i}\right]}f\left(s_{k-1}\right)\left[g\left(s_{k}\right)-g\left(s_{k-1}\right)\right]\right|.
\]
We estimate 
\begin{eqnarray*}
\lefteqn{\left|f\left(\nu_{i}\right)\left[g\left(a_{i}\right)-g\left(a_{i-1}\right)\right]-\sum_{k:s_{k-1},s_{k}\in\left[a_{i-1};a_{i}\right]}f\left(s_{k-1}\right)\left[g\left(s_{k}\right)-g\left(s_{k-1}\right)\right]\right|}\\
 & \leq & \left|f\left(\nu_{i}\right)\left[g\left(a_{i}\right)-g\left(a_{i-1}\right)\right]-f\left(a_{i-1}\right)\left[g\left(a_{i}\right)-g\left(a_{i-1}\right)\right]\right|\\
 &  & +\left|\sum_{k:s_{k-1},s_{k}\in\left[a_{i-1};a_{i}\right]}f\left(s_{k-1}\right)\left[g\left(s_{k}\right)-g\left(s_{k-1}\right)\right]-f\left(a_{i-1}\right)\left[g\left(a_{i}\right)-g\left(a_{i-1}\right)\right]\right|.
\end{eqnarray*}

Recall the definition of  $S.$ If there exists $N=0,1,2,\ldots$ such that $\eta_{N}=0$ or $\theta_{N}=0$ then $\TTV f{\left[a;b\right]}{}$ or $\TTV g{\left[a;b\right]}{}$ is finite, moreover, both functions $f$ and $g$ are bounded (since they are regulated), hence the integral $\int_a^b f\dd g$ exists. Thus we may and will assume that 
$\eta_{N}>0$ and $\theta_{N}>0$ for all $N=0,1,2,\ldots$

Choose $N=1,2,\ldots.$ By the assumption that $f$ and $g$ have
no common points of discontinuity, for sufficiently small $\mbox{mesh}\left(\pi\right),$
for $i=1,2,\ldots,l$ we have 
\begin{equation}
\sup_{a_{i-1}\leq s\leq a_{i}}\left|f\left(s\right)-f\left(a_{i-1}\right)\right|\leq\eta_{N-1}\label{eq:ind_f}
\end{equation}
or 
\begin{equation}
\sup_{a_{i-1}\leq s\leq a_{i}}\left|g\left(a_{i}\right)-g\left(s\right)\right|\leq\theta_{N-1}.\label{eq:ind_J}
\end{equation}
To see this, assume that for every $h>0,$ there exist $\left[a_{h};b_{h}\right]\subset\left[a;b\right]$
such that $b_{h}-a_{h}\leq h$ and $\sup_{x,y\in\left[a_{h};b_{h}\right]}\left|f\left(y\right)-f\left(y\right)\right|>\eta_{N-1}$
and $\sup_{x,y\in\left[a_{h};b_{h}\right]}\left|g\left(x\right)-g\left(y\right)\right|>\theta_{N-1}.$
We choose a convergent subsequence of the sequence $\left(a_{1/n}+b_{1/n}\right)/2,$
$n=1,2,.\ldots,$ and we see that the limit of this sequence is a
point of discontinuity for both $f$ and $g,$ which is a contradiction
with the assumption that $f$ and $g$ have no common points of discontinuity.\par Let
$I$ be the set of all indices $i=1,2,\ldots,l$ for which (\ref{eq:ind_f})
holds. Now, for $i\in I,$ set $\delta_{j-1}:=\eta_{N+j-1},$ $\varepsilon_{j}:=\theta_{N+j},$ $j=0,1,2,\ldots,$ and 
define 
\[
S_{i}:=\sum_{j=0}^{+\infty}2^{j}\eta_{j-1}\cdot\TTV g{\left[a_{i-1};a_{i}\right]}{\theta_{j}}+\sum_{j=0}^{+\infty}2^{j}\theta_{j}\cdot\TTV f{\left[a_{i-1};a_{i}\right]}{\eta_{j}}.
\]
By Lemma \ref{lema} we estimate 
\begin{eqnarray*}
 &  & \left|\sum_{k:s_{k-1},s_{k}\in\left[a_{i-1};a_{i}\right]}f\left(s_{k-1}\right)\left[g\left(s_{k}\right)-g\left(s_{k-1}\right)\right]-f\left(a_{i-1}\right)\left[g\left(a_{i}\right)-g\left(a_{i-1}\right)\right]\right|\\
 &  & \leq\sum_{j=0}^{+\infty}2^{j}\delta_{j-1}\cdot\TTV g{\left[a_{i-1};a_{i}\right]}{\varepsilon_{j}}+\sum_{j=0}^{+\infty}2^{j}\varepsilon_{j}\cdot\TTV f{\left[a_{i-1};a_{i}\right]}{\delta_{j}}\\
 &  & \leq\sum_{j=0}^{+\infty}2^{j}\eta_{N+j-1}\cdot\TTV g{\left[a_{i-1};a_{i}\right]}{\theta_{N+j}}+\sum_{j=0}^{+\infty}2^{j}\theta_{N+j}\cdot\TTV f{\left[a_{i-1};a_{i}\right]}{\eta_{N+j}}\\
 &  & \leq2^{-N}S_{i}.
\end{eqnarray*}
Similarly, 
\begin{eqnarray*}
\left|f\left(\nu_{i}\right)\left[g\left(a_{i}\right)-g\left(a_{i-1}\right)\right]-f\left(a_{i-1}\right)\left[g\left(a_{i}\right)-g\left(a_{i-1}\right)\right]\right| & \leq & 2^{-N}S_{i}.
\end{eqnarray*}
Hence 
\begin{eqnarray}
 &  & \left|f\left(\nu_{i}\right)\left[g\left(a_{i}\right)-g\left(a_{i-1}\right)\right]-\sum_{k:s_{k-1},s_{k}\in\left[a_{i-1};a_{i}\right]}f\left(s_{k-1}\right)\left[g\left(s_{k}\right)-g\left(s_{k-1}\right)\right]\right|\nonumber \\
 &  & \leq2^{1-N}S_{i}.\label{eq:ineq_i}
\end{eqnarray}
The truncated variation is a superadditive functional of the interval, from
which we have 
\[
\sum_{i\in I}\TTV g{\left[a_{i-1};a_{i}\right]}{\theta_{j}}\leq\TTV g{\left[a;b\right]}{\theta_{j}},
\]
\[
\sum_{i\in I}\TTV f{\left[a_{i-1};a_{i}\right]}{\eta_{j}}\leq\TTV f{\left[a;b\right]}{\eta_{j}}.
\]
By (\ref{eq:ineq_i}) and last two inequalities, summing over $i\in I$
we get the estimate 
\begin{eqnarray}
 &  & \left|\sum_{i\in I}\left\{ f\left(\nu_{i}\right)\left[g\left(a_{i}\right)-g\left(a_{i-1}\right)\right]-\sum_{k:s_{k-1},s_{k}\in\left[a_{i-1};a_{i}\right]}f\left(s_{k-1}\right)\left[g\left(s_{k}\right)-g\left(s_{k-1}\right)\right]\right\} \right|\nonumber \\
 &  & \leq2^{1-N}\sum_{i\in I}S_{i}\leq2^{1-N}S.\label{eq:estim_I-1}
\end{eqnarray}
\par Now, let $J$ be the set of all indices, for which (\ref{eq:ind_J})
holds. For $i=1,2,\ldots,l$ define 
\[
T_{i}:=\sum_{j=0}^{+\infty}2^{j}\theta_{j}\cdot\TTV f{\left[a_{i-1};a\right]}{\eta_{j}}+\sum_{j=0}^{+\infty}2^{j}\eta_{j}\cdot\TTV g{\left[a_{i-1};a_{i}\right]}{\theta_{j+1}}.
\]
For $i\in J,$ by the summation by parts and then by Lemma \ref{lema}
we get\par 
\begin{eqnarray*}
 &  & \left|f\left(a_{i}\right)\left[g\left(a_{i}\right)-g\left(a_{i-1}\right)\right]-\sum_{k:s_{k-1},s_{k}\in\left[a_{i-1};a_{i}\right]}f\left(s_{k-1}\right)\left[g\left(s_{k}\right)-g\left(s_{k-1}\right)\right]\right|\\
 &  & =\left|\sum_{k:s_{k-1},s_{k}\in\left[a_{i-1};a_{i}\right]}g\left(s_{k}\right)\left[f\left(s_{k}\right)-f\left(s_{k-1}\right)\right]-g\left(a_{i-1}\right)\left[f\left(a_{i}\right)-f\left(a_{i-1}\right)\right]\right|\\
 &  & \leq\sum_{j=0}^{+\infty}2^{j}\theta_{N+j-1}\cdot\TTV f{\left[a_{i-1};a_{i}\right]}{\eta_{N+j}}+\sum_{j=0}^{+\infty}2^{j}\eta_{N+j}\cdot\TTV g{\left[a_{i-1};a_{i}\right]}{\theta_{N+j}}.\\
 &  & \leq2^{1-N}T_{i}\leq2^{1-N}S_{i}.
\end{eqnarray*}
Similarly, by Lemma \ref{lema}, 
\begin{eqnarray*}
 &  & \left|f\left(a_{i}\right)\left[g\left(a_{i}\right)-g\left(a_{i-1}\right)\right]-f\left(\nu_{i}\right)\left[g\left(a_{i}\right)-g\left(a_{i-1}\right)\right]\right|\\
 &  & =\left|g\left(a_{i-1}\right)\left[f\left(\nu_{i}\right)-f\left(a_{i-1}\right)\right]+g\left(a_{i}\right)\left[f\left(a_{i}\right)-f\left(\nu_{i}\right)\right]-g\left(a_{i-1}\right)\left[f\left(a_{i}\right)-f\left(a_{i-1}\right)\right]\right|\\
 &  & \leq2^{1-N}S_{i}.
\end{eqnarray*}
From last two inequalities we get 
\begin{eqnarray*}
\lefteqn{\left|f\left(\nu_{i}\right)\left[g\left(a_{i}\right)-g\left(a_{i-1}\right)\right]-\sum_{k:s_{k-1},s_{k}\in\left[a_{i-1};a_{i}\right]}f\left(s_{k-1}\right)\left[g\left(s_{k}\right)-g\left(s_{k-1}\right)\right]\right|}\\
 & \leq & \left|f\left(a_{i}\right)\left[g\left(a_{i}\right)-g\left(a_{i-1}\right)\right]-f\left(\nu_{i}\right)\left[g\left(a_{i}\right)-g\left(a_{i-1}\right)\right]\right|\\
 &  & +\left|f\left(a_{i}\right)\left[g\left(a_{i}\right)-g\left(a_{i-1}\right)\right]-\sum_{k:s_{k-1},s_{k}\in\left[a_{i-1};a_{i}\right]}f\left(s_{k-1}\right)\left[g\left(s_{k}\right)-g\left(s_{k-1}\right)\right]\right|\\
 & \leq & 2^{2-N}S_{i}.
\end{eqnarray*}
Summing over $i\in J$ and using the superadditivity of the truncated
variation as a function of the interval, we get the estimate 
\begin{eqnarray}
 &  & \left|\sum_{i\in J}\left\{ f\left(\nu_{i}\right)\left[g\left(a_{i}\right)-g\left(a_{i-1}\right)\right]-\sum_{k:s_{k-1},s_{k}\in\left[a_{i-1};a_{i}\right]}f\left(s_{k-1}\right)\left[g\left(s_{k}\right)-g\left(s_{k-1}\right)\right]\right\} \right|\nonumber \\
 &  & \leq2^{2-N}\sum_{i\in J}S_{i}\leq2^{2-N}S.\label{eq:estim_J-1}
\end{eqnarray}
Finally, from (\ref{eq:estim_I-1}) and (\ref{eq:estim_J-1}) we get
\[
\left|f\left(\nu_{i}\right)\left[g\left(b\right)-g\left(a\right)\right]-\sum_{k=1}^{n}f\left(s_{k-1}\right)\left[g\left(s_{k}\right)-g\left(s_{k-1}\right)\right]\right|\leq6\cdot2^{-N}S.
\]
\par Similar estimate holds for 
\[
\left|\sum_{i=j}^{m}f\left(\xi_{j}\right)\left[g\left(b_{j}\right)-g\left(b_{j-1}\right)\right]-\sum_{k=1}^{n}f\left(s_{k-1}\right)\left[g\left(s_{k}\right)-g\left(s_{k-1}\right)\right]\right|,
\]
provided that $\mbox{mesh}\left(\rho\right)$ is sufficiently small.
Hence 
\[
\left|\sum_{i=1}^{l}f\left(\nu_{i}\right)\left[g\left(a_{i}\right)-g\left(a_{i-1}\right)\right]-\sum_{i=j}^{m}f\left(\xi_{j}\right)\left[g\left(b_{j}\right)-g\left(b_{j-1}\right)\right]\right|\leq12\cdot2^{-N}S.
\]
provided that $\mbox{mesh}\left(\pi\right)$ and $\mbox{mesh}\left(\rho\right)$
are sufficiently small. Since $N$ may be arbitrary large, we get
the convergence of the approximating sums to an universal limit, which
is the Riemann-Stieltjes integral. The estimate (\ref{eq:estimate_integral})
follows directly from the proved convergence of approximating sums
to the Riemann-Stieltjes integral and Lemma \ref{lema}. 
\end{dwd} 
Using Remark \ref{symmetry} and reasoning similarly as in the proof of Theorem \ref{main}, we get the symmetric result.
\begin{theo} \label{main1} Let $f,g:\left[a;b\right]\rightarrow\R$
be two regulated functions which have no common points of discontinuity.
Let $\eta_{0}\geq\eta_{1}\geq\ldots$ and $\theta_{0}\geq\theta_{1}\geq\ldots$
be two sequences of non-negative numbers, such that $\eta_{k}\downarrow0,$
$\theta_{k}\downarrow0$ as $k\rightarrow+\infty.$ Define $\theta_{-1}:=\sup_{a\leq t\leq b}\left|g\left(b\right)-g\left(t\right)\right|$
and 
\[
\tilde{S}:=\sum_{k=0}^{+\infty}2^{k}\theta_{k-1}\cdot\TTV f{\left[a;b\right]}{\eta_{k}}+\sum_{k=0}^{\infty}2^{k}\eta_{k}\cdot\TTV g{\left[a;b\right]}{\theta_{k}}.
\]
If $\tilde{S}<+\infty$ then the Riemann-Stieltjes integral 
$
\int_{a}^{b}f\mathrm{d}g
$
exists and one has the following estimate 
\begin{equation}
\left|\int_{a}^{b}f\mathrm{d}g-f\left(a\right)\left[g\left(b\right)-g\left(a\right)\right]\right|\leq \tilde{S}.\label{eq:estimate_integral1}
\end{equation}
\end{theo}
From Theorem \ref{main}, Theorem \ref{main1} and Remark \ref{n=1} we also  have
\begin{coro} \label{minST} Let $f,g:\left[a;b\right]\rightarrow\R$
be two regulated functions which have no common points of discontinuity, $\xi \in [a;b]$ and $S$ and $\tilde{S}$ be as in Theorem \ref{main} and Theorem \ref{main1} respectively. If $S<+\ns$ or $\tilde{S}< +\ns$ then the Riemann-Stieltjes integral $\int_{a}^{b}f\mathrm{d}g$ exists and one has the following estimate 
\[
\left|\int_{a}^{b}f\mathrm{d}g-f\left(\xi\right)\left[g\left(b\right)-g\left(a\right)\right]\right|\leq 
2 \min\{S,\tilde{S} \}.
\]
\end{coro}

\subsection{Young's Theorem and the Lo\'{e}ve-Young inequality}

Let for $p>0,$ ${\cal V}^{p}\rbr{[a;b]}$ denote the family of functions $f:[a;b]\ra \R$ with finite $p-$variation. Note that if $f \in {\cal V}^{p}\rbr{[a;b]}$ then $f$ is regulated. The additional relation we will use, is the following one: if $f\in{\cal V}^{p}\rbr{[a;b]}$ for some
$p\geq 1,$ then for every $\delta>0,$ 
\begin{equation}
\TTV f{\left[a;b\right]}{\delta}\leq V^{p}\left(f,\left[a;b\right]\right)\delta^{1-p}.\label{eq:p_variation}
\end{equation}
As far as we know, the first result of this kind, namely, $\TTV f{\left[a;b\right]}{\delta}\leq C_{f}\delta^{1-p}$
for a continuous function $f\in{\cal V}^{p}\rbr{[a;b]}$ and some
constant $C_{f}<+\infty$ depending on $f,$ 
was proven in \cite[Sect. 6]{TronelVladimirov:2000}. In \cite{TronelVladimirov:2000},
$\TTV f{\left[a;b\right]}{\varepsilon}$ is called $\varepsilon-$variation
and is denoted by $V_{f}(\varepsilon).$ However, being equipped with formula (\ref{TV_def1}) we see that relation (\ref{eq:p_variation})
follows immediately from the inequality: for any $a\leq s<t\leq b,$
\[
\max\left\{ \left|f\left(t\right)-f\left(s\right)\right|-\delta,0\right\} \leq\frac{\left|f\left(t\right)-f\left(s\right)\right|^{p}}{\delta^{p-1}},
\]
which is an obvious consequence of the estimate:
\[
\delta^{p-1}\max\left\{ \left|x\right|-\delta,0\right\} \leq\begin{cases}
0 & \mbox{if }\delta\geq\left|x\right|\\
\left|x\right|^{p-1}\max\left\{ \left|x\right|-\delta,0\right\}  & \mbox{if }0<\delta<\left|x\right|
\end{cases}\leq\left|x\right|^{p}
\]
for any $\delta >0$ and any real $x.$ 

Let us denote 
\begin{equation} \label{p_var_norm_def}
\left\Vert f\right\Vert_{p-\text{var},\left[a;b\right]}:=\rbr{V^{p}\left(f,\left[a;b\right]\right)}^{1/p}
\end{equation}
and recal that $\left\Vert f\right\Vert_{\text{osc},\left[a;b\right]}=\sup_{a\leq s<t\leq b}
\left|f\left(t\right)-f\left(s\right)\right|.$
Now we are ready to state a Corollary stemming from Theorem \ref{main}, which was one of the main results of \cite{Young:1936}. The second part of this Corollary is an improved version of the Lo\'{e}ve-Young inequality. 
\begin{coro} \label{corol_Young} Let $f,g:\left[a;b\right]\rightarrow\R$
be two functions with no common points of discontinuity.
If $f\in{\cal V}^{p}\rbr{[a;b]}$ and $g\in{\cal V}^{q}\rbr{[a;b]},$ where $p>1,$ $q>1,$
$p^{-1}+q^{-1}>1,$ then the Riemann Stieltjes $\int_{a}^{b}f\mathrm{d}g$
exists. Moreover, there exist a constant $C_{p,q},$ depending on $p$ and $q$ only, such that  
\begin{align*}
\left|\int_{a}^{b}f\mathrm{d}g-f\left(a\right)\left[g\left(b\right)-g\left(a\right)\right]\right| 
  \leq C_{p,q}\left\Vert f\right\Vert_{p-\emph{var},\left[a;b\right]}^{p-p/q}\left\Vert f\right\Vert _{\emph{osc},\left[a;b\right]}^{1+p/q-p} \left\Vert g\right\Vert_{q-\emph{var},\left[a;b\right]}.
\end{align*}
\end{coro} 
\begin{dwd} By Theorem \ref{main} it is enough to prove that for
some positive sequences $\eta_{0}\geq \eta_1 \geq \ldots$
and $\theta_{0}\geq\theta_{1}\geq\ldots,$ such that $\eta_{k}\downarrow0,$
$\theta_{k}\downarrow0$ as $k\rightarrow+\infty$ and  $\eta_{-1}=\sup_{a\leq t\leq b}\left|f\left(t\right)-f\left(a\right)\right|$ one has 
\begin{align*}
S: & =\sum_{k=0}^{+\infty}2^{k}\eta_{k-1}\cdot\TTV g{\left[a;b\right]}{\theta_{k}}+\sum_{k=0}^{+\infty}2^{k}\theta_{k}\cdot\TTV f{\left[a;b\right]}{\eta_{k}},\\
 & \leq C_{p,q}\left\Vert f\right\Vert_{p-\text{var},\left[a;b\right]}^{p-p/q}\left\Vert f\right\Vert _{\text{osc},\left[a;b\right]}^{1+p/q-p} \left\Vert g\right\Vert_{q-\text{var},\left[a;b\right]}.
\end{align*}
The proof will follow from the proper choice of $\left(\eta_{k}\right)$
and $\left(\theta_{k}\right).$ Since $p^{-1}+q^{-1}>1,$ we have
$\left(q-1\right)\left(p-1\right)<1.$ We choose \[\alpha \in \rbr{\sqrt{(q-1)(p-1)};1}, \quad 
\beta=\sup_{a\leq t\leq b}\left|f\left(t\right)-f\left(a\right)\right|, \quad  \gamma>0
\] and for $k=0,1,\ldots,$ define 
\[
\eta_{k-1}=\beta \cdot  2^{-\left(\alpha^{2}/\left[\left(q-1\right)\left(p-1\right)\right]\right)^{k}+1},
\]
\[
\theta_{k}= \gamma \cdot 2^{-\left(\alpha^{2}/\left[\left(q-1\right)\left(p-1\right)\right]\right)^{k}\alpha/\left(q-1\right)}.
\]
By (\ref{eq:p_variation}) we estimate 
\begin{align*}
\eta_{k-1}\cdot\TTV g{\left[a;b\right]}{\theta_{k}}\leq & \beta \cdot 2^{-\left(\alpha^{2}/\left[\left(q-1\right)\left(p-1\right)\right]\right)^{k}+1}\\
 & \times V^{q}\left(g,\left[a;b\right]\right)\left(\gamma \cdot 2^{-\left(\alpha^{2}/\left[\left(q-1\right)\left(p-1\right)\right]\right)^{k}\alpha/\left(q-1\right)}
\right)^{1-q}\\
= & 2^{-\left(1-\alpha\right)\left(\alpha^{2}/\left[\left(q-1\right)\left(p-1\right)\right]\right)^{k}+1}V^{q}\left(g,\left[a;b\right]\right)\beta\gamma^{1-q},
\end{align*}
and similarly 
\begin{align*}
\theta_{k}\cdot\TTV f{\left[a;b\right]}{\eta_{k}}\leq & \gamma \cdot 2^{-\left(\alpha^{2}/\left[\left(q-1\right)\left(p-1\right)\right]\right)^{k}\alpha/\left(q-1\right)}
\\
 & \times V^{p}\left(f,\left[a;b\right]\right)\left(\beta \cdot 2^{-\left(\alpha^{2}/\left[\left(q-1\right)\left(p-1\right)\right]\right)^{k+1}+1}
\right)^{1-p}\\
= & 2^{-\left(1-\alpha\right)\left(\alpha^{2}/\left[\left(q-1\right)\left(p-1\right)\right]\right)^{k}\alpha/\left(q-1\right)+1-p}V^{p}\left(f,\left[a;b\right]\right)
\beta^{1-p}\gamma.
\end{align*}
Hence 
\begin{eqnarray*}
S & = & \sum_{k=0}^{+\infty}2^{k}\eta_{k-1}\cdot\TTV g{\left[a;b\right]}{\theta_{k}}+\sum_{k=0}^{+\infty}2^{k}\theta_{k}\cdot\TTV f{\left[a;b\right]}{\eta_{k}}\\
 & \leq & \left(\sum_{k=0}^{+\infty}2^{k}2^{-\left(1-\alpha\right)\left(\alpha^{2}/\left[\left(q-1\right)\left(p-1\right)\right]\right)^{k}+1}\right)V^{q}\left(g,\left[a;b\right]\right)\beta\gamma^{1-q}\\
 &  & +\left(\sum_{k=0}^{+\infty}2^{k}2^{-\left(1-\alpha\right)\left(\alpha^{2}/\left[\left(q-1\right)\left(p-1\right)\right]\right)^{k}\alpha/\left(q-1\right)+1-p}\right)V^{p}\left(f,\left[a;b\right]\right)
\beta^{1-p}\gamma.
\end{eqnarray*}
Since $\alpha<1$ and $\alpha^{2}/\left[\left(q-1\right)\left(p-1\right)\right]>1,$
we easily infer that $S<+\infty,$ from which we get that the integral
$\int_{a}^{b}f\mathrm{d}g$ exists. 

Moreover,
denoting 
\begin{align*}
C_{p,q} & = \max\left\{ \sum_{k=0}^{+\infty}2^{k+2-\left(1-\alpha\right)\left(\alpha^{2}/\left[\left(q-1\right)\left(p-1\right)\right]\right)^{k}}, \right. \\
& \quad \quad \quad \quad \left.\sum_{k=0}^{+\infty}2^{k+2-\left(1-\alpha\right)\left(\alpha^{2}/\left[\left(q-1\right)\left(p-1\right)\right]\right)^{k}\alpha/\left(q-1\right)-p}\right\} 
\end{align*}
we get 
\begin{align*}
S & \leq\frac{1}{2}C_{p,q}\left(V^{q}\left(g,\left[a;b\right]\right)\beta\gamma^{1-q}+V^{p}\left(f,\left[a;b\right]\right)
\beta^{1-p}\gamma\right).
\end{align*}
Setting in this expression 
$
\gamma=\left(V^{q}\left(g,\left[a;b\right]\right)\right/V^{p}\left(f,\left[a;b\right]\right))^{1/q}
\beta^{p/q}
$
we obtain 
\begin{align*}
S & \leq C_{p,q}\left(V^{q}\left(g,\left[a;b\right]\right)\right)^{1/q}\left(V^{p}\left(f,\left[a;b\right]\right)\right)^{1-1/q}
\beta^{1+p/q-p}
\\ 
& \leq C_{p,q}\left\Vert g\right\Vert_{q-\text{var},\left[a;b\right]}\left\Vert f\right\Vert_{p-\text{var},\left[a;b\right]}^{p-p/q}\left\Vert f\right\Vert _{\text{osc},\left[a;b\right]}^{1+p/q-p}.
\end{align*}
\end{dwd} 
\begin{rem} \label{rem_Young}
Let $f,$ $g,$ $p,$ $q$ and $C_{p,q}$ be the same as in Corollary \ref{corol_Young}. Using Theorem \ref{main1} instead of Theorem \ref{main}, we get the following, similar estimate
\begin{align*}
\left|\int_{a}^{b}f\mathrm{d}g-f\left(a\right)\left[g\left(b\right)-g\left(a\right)\right]\right| 
  \leq C_{p,q}\left\Vert f\right\Vert_{p-\emph{var},\left[a;b\right]}\left\Vert g\right\Vert_{q-\emph{var},\left[a;b\right]}^{q-q/p}\left\Vert g\right\Vert _{\emph{osc},\left[a;b\right]}^{1+q/p-q}.
\end{align*}
From Corollary \ref{minST} and the obtained estimates, we also have that for any $\xi \in [a;b]$ 
\begin{align*}
& \left|\int_{a}^{b}f\mathrm{d}g-f\left(\xi\right)\left[g\left(b\right)-g\left(a\right)\right]\right| 
 \\
& \leq 2 C_{p,q}
\left\Vert f\right\Vert_{p-\emph{var},\left[a;b\right]} \left\Vert g\right\Vert_{q-\emph{var},\left[a;b\right]}
\min \cbr{\frac{\left\Vert f\right\Vert _{\emph{osc},\left[a;b\right]}^{1+p/q-p}}{\left\Vert f\right\Vert_{p-\emph{var},\left[a;b\right]}^{1+p/q-p} 
}, 
\frac{\left\Vert g\right\Vert _{\emph{osc},\left[a;b\right]}^{1+q/p-q}}{{\left\Vert g\right\Vert_{q-\emph{var},\left[a;b\right]}^{1+q/p-q}}} 
}.
\end{align*}
\end{rem}
\begin{rem} \label{rem_int_V^p}
From Corollary \ref{corol_Young}, reasoning in the similar way as in \cite[p. 456]{Lyons:1994}, we get the following important estimate of the $q$-variation of the function $t \mapsto \int_a^{t} f \dd g$
\begin{eqnarray*}
\Varnormthm {\int_a^{\cdot} f \dd g}{q} & \leq & \rbr{ C_{p,q}\left\Vert f\right\Vert_{p-\emph{var},\left[a;b\right]}^{p-p/q}\left\Vert f\right\Vert _{\emph{osc},\left[a;b\right]}^{1+p/q-p} + \left\Vert f\right\Vert_{\ns, [a;b]} }\left\Vert g\right\Vert_{q-\emph{var},\left[a;b\right]}, \\
& \leq & \rbr{ C_{p,q}\left\Vert f\right\Vert_{p-\emph{var},\left[a;b\right]} + \left\Vert f\right\Vert_{\ns, [a;b]} }\left\Vert g\right\Vert_{q-\emph{var},\left[a;b\right]},
\end{eqnarray*}
where $f,$ $g,$ $p,$ $q$ and $C_{p,q}$ are the same as in Corollary \ref{corol_Young}.
\end{rem}

\section{Integral equations driven by moderately irregular signals}
Let $p \in (1;2).$ The preceding section provides us with tools to solve integral equations of the following form
\begin{equation} \label{eq:integral}
y(t)=y_0 + \int_a^t F(y(s)) \dd x(s),
\end{equation}
where $x$ is a continuous function from the space ${\cal V}^p\rbr{[a;b]}$  and $F:\R \ra \R$ is $\alpha$-Lipschitz. 
The  functional $\Varfullnorm{\cdot}{p} : {\cal V}^p\rbr{[a;b]} \ra [0;+\ns)$ defined as
\[
\Varfullnorm{f}{p} := |f(a)|+\Varnorm{f}{p}
\] is a norm and the space ${\cal V}^p\rbr{[a;b]}$ equipped with this norm is a Banach space. 
For our purposes it will be enough to work with the following definition
of locally or globally $\alpha$-Lipschitz function when $\alpha\in\left(0;1\right]$.
For $x=\left(x_{1},\ldots,x_{n}\right)\in\R^{n}$ we denote $\left\Vert x\right\Vert =\max_{i=1,\ldots,n}\left|x_{i}\right|.$
\begin{defi}
Let $F:\R^{n}\ra\R$ and $\alpha\in\left(0;1\right].$
For any $R>0$ we define its \emph{local $\alpha$-Lipschitz parameter
$K_{F}^{\alpha}\left(R\right)$ as} 
\[
K_{F}^{\left(\alpha\right)}\left(R\right):=\sup\left\{ \frac{\left|F\left(y\right)-F\left(x\right)\right|}{\left\Vert y-x\right\Vert ^{\alpha}}:x,y\in\R^{n},x\neq y,\left\Vert x\right\Vert \leq R,\left\Vert y\right\Vert \leq R\right\} 
\]
and its \emph{global $\alpha$-Lipschitz parameter $K_{F}^{\left(\alpha\right)}$
as }$K_{F}^{\left(\alpha\right)}:=\lim_{R\ra+\ns}K_{F}\left(R\right)<+\ns.$
The function $F$ will be called \emph{locally $\alpha$-Lipschitz}
if for every $R>0,$ $K_{F}^{\left(\alpha\right)}\left(R\right)<+\ns$
and it will be called \emph{globally $\alpha$-Lipschitz}
if $K_{F}^{\left(\alpha\right)}<+\ns.$ 
\end{defi}
In the case, when there is no ambiguity what is the value of the parameter
$\alpha$ and what is the function $F,$ we will write $K_{F}\left(R\right),$
$K_{F},$ $K\left(R\right)$  or even $K.$

First we will consider the case $p-1 < \alpha < 1.$ In this case we have the existence but no uniqueness result. We will obtain a stronger result than similar results \cite[Lemma, p. 459]{Lyons:1994} or \cite[Theorem 1.20]{LyonsCaruana:2007}. Namely, we will prove that there exists a solution to (\ref{eq:integral}) {\bf which is an element of the space} ${\cal V}^p\rbr{[a;b]},$ not only an element of the space ${\cal V}^{q}\rbr{[a;b]}$ for arbitrary chosen $q>p.$ This will be possible with the use of Remark \ref{rem_int_V^p}.
\begin{prop} \label{prop:int_equation}
Let $p\in(1;2),$ $y_0\in \R,$ $x$ be a continuous function from the space ${\cal V}^p\rbr{[a;b]}$ and $F:\R\ra \R$ be globally $\alpha$-Lipschitz where $p-1 < \alpha < 1.$ Equation (\ref{eq:integral}) admits a solution $y,$ which is an element of ${\cal V}^p\rbr{[a;b]}.$ Moreover, $\Varfullnormthm{y}{p}\leq R,$ where $R>0$ satisfies the equality 
\[
R = \rbr{C_{p/\alpha,p} +2}K^{(\alpha)}_F\Varnormthm{x}{p} R^{\alpha} + \left| y_0 \right| +\left| F(0)\right|\Varnormthm{x}{p},
\]
with $C_{p/\alpha,p}$ being the same as in Corollary \ref{corol_Young} and Remark \ref{rem_int_V^p}.
\end{prop}
Now we proceed to the proof of Proposition \ref{prop:int_equation}. We will proceed in a standard way, but with the more accurate estimate of Remark \ref{rem_int_V^p} we will be able to obtain the finiteness of $\Varfullnorm {\cdot}{p}$ norm of the solution.
\begin{dwd}
Let $f \in {\cal V}^{p}\rbr{[a;b]}.$ By \cite[Lemma 1.18]{LyonsCaruana:2007}, $F(f(\cdot)) \in {\cal V}^{p/\alpha}\rbr{[a;b]}$ and since $\alpha/p+1/p>1,$ we may apply Remark \ref{rem_int_V^p} and define the operator $T: {\cal V}^{p}\rbr{[a;b]} \ra {\cal V}^{p}\rbr{[a;b]},$ 
\[
Tf := y_0 + \int_a^{\cdot} F(f(t))\dd x(t).
\] 
Denote $K=K^{(\alpha)}_F.$ Using Remark \ref{rem_int_V^p} and \cite[Lemma 1.18]{LyonsCaruana:2007} we estimate 
\begin{align}
& \Varfullnorm{Tf}{p} = \Varfullnorm {y_0 + \int_a^{\cdot} F(f(t))\dd x(t)}{p} \nonumber \\ 
& \leq \left| y_0 \right| + \Varnorm {\int_a^{\cdot} \sbr{ F(f(t)) - F(f(a)) }\dd x(t)}{p} + \left| F(f(a))\right|\Varnorm {x}{p} \nonumber \\
& \leq \left| y_0 \right| + \rbr{ C_{p/\alpha,p} \Varnorm {F(f(\cdot))}{p/\alpha}+ \left\Vert F(f(\cdot)) - F(f(a))\right\Vert_{\ns, [a;b]} + \left| F(f(a))\right| } \Varnorm {x}{p} \nonumber \\ 
& \leq \left| y_0 \right| + \rbr{ \rbr{C_{p/\alpha,p}+1} \Varnorm {F(f(\cdot))}{p/\alpha}+ \left| F(f(a))\right| } \Varnorm {x}{p} \nonumber \\ 
& \leq \left| y_0 \right| + \rbr{ \rbr{C_{p/\alpha,p}+1} K \Varnorm {f}{p}^{\alpha} + \left| F(f(a))\right| }\Varnorm {x}{p}. \label{ineq:TVnorm_est}
\end{align} 
By the Lipschitz property, 
\[
\left| F(f(a))\right| \leq K\left| f(a)\right|^{\alpha} + \left| F(0)\right| \leq K\Varfullnorm {f}{p}^{\alpha} + \left| F(0)\right|.
\]
Denoting 
\begin{equation} \label{defAB} 
A = \rbr{C_{p/\alpha,p}+2}K\Varnorm {x}{p} \text{ and } B= \left| y_0 \right| +\left| F(0)\right|\Varnorm {x}{p},
\end{equation}
from (\ref{ineq:TVnorm_est}) we get 
\begin{equation} \label{ineq:stabil}
\Varfullnorm {Tf}{p} \leq A \Varfullnorm {f}{p}^{\alpha} + B.
\end{equation}

For $\alpha <1$ let $R$ be the least positive solution of the inequality $R \geq  A \cdot R^\alpha + B$ (i.e. $R = A \cdot R^\alpha + B$). From (\ref{ineq:stabil}) we have that the operator $T$ maps the closed ball ${\cal B}(R) = \cbr{f\in {\cal V}^{p}\rbr{[a;b]} : \Varfullnorm{f}{p} \leq R}$ to itself. 

Now, for $f,g \in {\cal V}^{p}\rbr{[a;b]}$ we are going to investigate the difference $Tf-Tg.$ Again, using Remark \ref{rem_int_V^p} and the Lipschitz property we estimate
\begin{align}
& \Varfullnorm{Tf - Tg}{p} = \Varnorm {\int_a^{\cdot} \sbr{ F(f(t)) -F(g(t))}\dd x(t)}{p}  \nonumber \\ 
& \leq \Varnorm {\int_a^{\cdot} \sbr{ F(f(t)) -F(g(t)) - \cbr{F(f(a)) -F(g(a))}}\dd x(t)}{p} \nonumber \\ & \quad +  \left| F(f(a)) -F(g(a)) \right| \Varnorm {x}{p}  \nonumber \\
& \leq \rbr{C_{p/\alpha,p} +1} \Varnorm {F(f(\cdot))-F(g(\cdot))}{p/\alpha}^{(p-1)/\alpha} \Oscnorm {F(f(\cdot))-F(g(\cdot))}^{(\alpha+1-p)/\alpha} \Varnorm {x}{p} \nonumber \\ & \quad + \left| F(f(a)) -F(g(a)) \right| \Varnorm {x}{p} \nonumber \\
&  \leq  \rbr{C_{p/\alpha,p} +1} K \rbr{\Varnorm {f}{p}^{\alpha}+\Varnorm {g}{p}^{\alpha}}^{(p-1)/\alpha} \Oscnorm {f-g}^{\alpha+1-p} \Varnorm {x}{p} \nonumber \\
& \quad +   K \left|f(a) -g(a) \right|^{\alpha}  \Varnorm {x}{p}. \label{lasttt}
\end{align} 
From (\ref{lasttt}) we see that $T$ is continuous. Moreover, 
from the first inequality in Remark \ref{rem_Young} and the continuity of $x$ we get that functions belonging to the image $T({\cal B}(R))$ are equicontinuous. Let ${\cal U}$ be the closure of the convex hull of $T({\cal B}(R))$ (in the topology induced by the norm $\Varfullnorm{\cdot}{p}$). It is easy to see that functions belonging to ${\cal U}$ are also equicontinuous. Moreover, ${\cal U} \subset {\cal B}(R)$ (since $T({\cal B}(R)) \subset {\cal B}(R)$ and ${\cal B}(R)$ is convex) and $T({\cal U}) \subset {\cal U}$ (since $T({\cal U}) \subset T({\cal B}(R)) $). Now, let ${\cal V} = T({\cal U}).$ From the equicontinuity of ${\cal U},$ Arzela-Ascoli Theorem and (\ref{lasttt}) we see that the set ${\cal V}$ is compact in the topology induced by the norm $\Varfullnorm{\cdot}{p}.$  Thus, by the fixed-point Theorem of Schauder, we get that there exists a point $y \in {\cal U}$ such that $Ty=y.$
\end{dwd}

Now we will consider the case $\alpha=1.$ 
\begin{fact}
Let $p\in(1;2),$ $y_0\in \R,$ $x$ be a continuous function from the space ${\cal V}^p\rbr{[a;b]}$ and $F:\R\ra \R$ be globally $1$-Lipschitz. Equation (\ref{eq:integral}) admits a solution $y,$ which is an element of ${\cal V}^p\rbr{[a;b]}.$ 
\end{fact}
\begin{dwd} 
To prove the assertion we may proceed in a similar way as in the proof of Proposition \ref{prop:int_equation}. The only thing we need is to assure that the inequality $R\geq A\cdot R + B,$ where $A$ and $B$ are defined in display (\ref{defAB}), holds for sufficiently large $R.$ This may be achieved by splitting the interval $[a;b]$ into small intervals, such that $A < 1$ on each of these intervals, and then solving equation (\ref{eq:integral}) on each of these intervals with the initial condition being equal the terminal value of the solution on the preceding interval. This is possible since for any $\varepsilon >0$ there exists $\delta >0$ such that for any $[c;c+\delta] \subset [a;b],$ ${\left\Vert {x} \right\Vert_{{p}-\text{var},\left[c;c+\delta\right]}} \leq \varepsilon,$ which is the consequence of the fact that the function $[a;b] \ni t \mapsto \left\Vert {x} \right\Vert_{{p}-\text{var},\left[a;t\right]}$ is continuous and we have  $\left\Vert {x} \right\Vert_{{p}-\text{var},\left[a;c+\delta\right]} \geq \left\Vert {x} \right\Vert_{{p}-\text{var},\left[a;c\right]} + \left\Vert {x} \right\Vert_{{p}-\text{var},\left[c;c+\delta\right]}$.  
\end{dwd}

\textbf{Acknowledgments} I would like to thank Raouf Ghomrasni
for drawing my attention to paper \cite{TronelVladimirov:2000} and Professor Terry Lyons for drawing my attention to his paper \cite{Lyons:1994}.

\bigskip

\textbf{Competing interests} The author declares that there are no competing  interests regarding the publication of this paper.


\begin{thebibliography}{10}

\bibitem{Banas_et_al:2013}
Appel, J, Bana\'{s},  J and Merentes, DNJ:
\newblock {\em Bounded variation and around}.
\newblock Series in
Nonlinear Analysis and Applications 17.
\newblock De Gruyter, Berlin (2013)

\bibitem{NorvaisaConcrete:2010}
Dudley, RM, Norvai\v{s}a, R:
\newblock {\em Concrete Functional Calculus}.
\newblock Springer Monographs in Mathematics. Springer, New York (2010)

\bibitem{Dyackov:1988}
D'ya\v{c}kov, AM: 
\newblock Conditions for the existence of {S}tieltjes integral of functions of
  bounded generalized variation.
\newblock {\em Anal. Math. (Budapest)} {\bf 14}, 295--313 (1988)

\bibitem{Friz:2010fk}
Friz, PK, Victoir, NB:
\newblock {\em Multidimensional stochastic processes as rough paths}, volume
  120 of {\em Cambridge Studies in Advanced Mathematics. Theory and
  applications}.
\newblock Cambridge University Press, Cambridge (2010)

\bibitem{LochowskiGhomrasniMMAS:2014}
{\L}ochowski, RM, Ghomrasni, R:
\newblock The play operator, the truncated variation and the generalisation of
  the {J}ordan decomposition.
\newblock {\em Math. Methods Appl. Sci.} {\bf 38}(3), 403-419 (2015)

\bibitem{LochowskiColloquium:2013}
{\L}ochowski, RM:
\newblock On the generalisation of the {Hahn-Jordan} decomposition for real
  c{\`a}dl{\`a}g functions.
\newblock {\em Colloq. Math.} {\bf 132}(1), 121--138 (2013)

\bibitem{Lyons:1994}
Lyons, TJ:
\newblock Differential equations driven by rough signals ({I}): an extension of
  an inequality of {L.} {C.} {Y}oung.
\newblock {\em Math. Res. Lett.} {\bf 1}, 451--464 (1994)

\bibitem{LyonsCaruana:2007}
Lyons, TJ, Caruana, M, L\'{e}vy, T:
\newblock {\em Differential Equations Driven By Rough Paths}, volume 1902 of
  {\em Lecture Notes in Mathematics}.
\newblock Springer Verlag, Berlin, Heidelberg (2007)

\bibitem{TronelVladimirov:2000}
Tronel, G, Vladimirov, AA:
\newblock On {BV}-type hysteresis operators.
\newblock {\em Nonlinear Anal. Ser. A: Theory Methods.} {\bf 39}(1), 79--98 (2000)

\bibitem{Young:1936}
Young, LC:
\newblock An inequality of the {H}\"{o}lder type, connected with {S}tieltjes
  integration.
\newblock {\em Acta Math. (Sweden)} {\bf 67}(1), 251--282 (1936)

\bibitem{Young:1938}
Young, LC:
\newblock General inequalities for {S}tieltjes integrals and the convergence of
  {F}ourier series.
\newblock {\em Math. Ann.} {\bf 115}, 581--612 (1938)

\end{thebibliography}
\end{document}